\documentclass[12pt]{amsart}

\usepackage[margin=1.3in]{geometry}
\usepackage[UKenglish]{babel} 
\usepackage[T1]{fontenc}
\usepackage[ansinew]{inputenc}
\usepackage{graphicx} 
\usepackage{lmodern} 
\usepackage{enumerate}
\usepackage{amsmath}
\usepackage{amsthm}
\usepackage{amsfonts}
\usepackage{amssymb}
\usepackage{amscd}
\usepackage{xcolor}
\usepackage{cleveref}
\usepackage[numbers,sort]{natbib}
\usepackage[colorinlistoftodos]{todonotes}
\usepackage{mathtools}

\newcommand{\newword}[1]{\emph{#1}}

\newcommand{\ldb}{\{\!\!\{}
\newcommand{\rdb}{\}\!\!\}}
\newcommand{\AAA}{\mathcal{A}}

\newcommand{\BB}{\mathcal{B}}
\newcommand{\CF}{\mathcal{CF}}

\newcommand{\MM}{\mathcal{M}}

\newcommand{\TT}{\mathbb{T}}

\newcommand{\cork}{\operatorname{cork}}

\newcommand{\rk}{\operatorname{rk}}

\newcommand{\stiefel}{\operatorname{St}} 
\newcommand{\dapx}{\mathcal{DA}} 
\newcommand{\dmat}{\mathcal{DM}} 
\newcommand{\f}[2]{#1_{#2}} 
\newcommand{\vertex}[2]{v^{#1}_{#2}} 
\newcommand{\Be}[1]{\mathcal{L}(#1)} 
\newcommand{\underm}[1]{\underbracket[1pt][-1pt]{#1}} 

\newcommand{\supp}{\operatorname{supp}}
\newcommand{\Moebius}{\operatorname{M\mbox{\"o}b}} 
\newcommand{\mult}{\operatorname{mult}} 

\newcommand{\nE}{E} 
\newcommand{\nnE}{E} 
\newcommand{\nEplusone}{E\sqcup\{*\}}
\newcommand{\nnEplusone}{E\sqcup\{*\}}
\newcommand{\nplusone}{*}

\theoremstyle{plain}
\newtheorem{proposition}{Proposition}[section]
\newtheorem{question}[proposition]{Question}
\newtheorem{theorem}[proposition]{Theorem}
\newtheorem{lemma}[proposition]{Lemma}

\theoremstyle{definition}
\newtheorem{definition}[proposition]{Definition}
\newtheorem{example}[proposition]{Example}

\theoremstyle{remark}
\newtheorem{remark}[proposition]{Remark}

\numberwithin{equation}{section}

\begin{document}

\pagestyle{plain}

\title{Extensions of transversal valuated matroids}
\author{Alex Fink, Jorge Alberto Olarte}
\date{}

\begin{abstract}
Following up on our previous work, we study single-element extensions of transversal valuated matroids.
We show that tropical presentations of valuated matroids with a minimal set of finite entries
enjoy counterparts of the properties proved by Bonin and de Mier of minimal non-valuated transversal presentations.
\end{abstract}

\maketitle

\section{Introduction}
Transversal matroids are a foundational class of matroids derived from bipartite graphs;
the graph is said to be a \emph{presentation} of the matroid.
Brualdi and Dinolt \cite{BrualdiDinolt}
studied the set of all presentations of a fixed transversal matroid $M$.
If we consider only presentations with the fewest vertices, 
where the part of the vertex set not labelled by the ground set of~$M$ 
has cardinality $\rk(M)$,
they showed this set has significant structure under the subgraph containment order.
Notably, in any \emph{minimal} presentation of~$M$, a vertex not in the ground set
has as its neighborhood a cocircuit of~$M$.
Starting with any presentation of~$M$,
and inserting into each of these neighbourhoods $N$ all coloops of the deletion $M\setminus N$,
gives the maximal presentation of~$M$, 
which Bondy~\cite{Bondy} and Mason~\cite{Mason-thesis} had earlier shown is unique.

Crapo \cite{Crapo} gave a classification of the single-element extensions of a general matroid,
but there remain many open problems about their behaviour.
A recent paper of Bonin and de Mier \cite{BdM} gives a wealth of structural information
about single-element extensions of a transversal matroid $M$ that are again transversal,
and the relationship of their presentations to presentations of~$M$.
Adding a vertex to the ground set side of a presentation of~$M$
produces a presentation of an extension of~$M$.
Potentially multiple extensions can be produced this way, 
depending on the neighbourhood chosen for the new vertex:
in fact all these extensions are distinct if and only if the starting presentation was minimal \cite[Theorem 3.1]{BdM}.
Moreover all transversal extensions of~$M$, aside from the single one that adds a coloop, 
arise from a minimal presentation of~$M$ this way \cite[Theorem 3.3]{BdM}.

The goal of the present work is to show comparable results for \emph{valuated matroids}. 
Valuated matroids, introduced by Dress and Wenzel \cite{DressWenzel},
are matroids enriched by the attachment of real-valued \emph{valuations},
to bases, elements of circuits, or in other cryptomorphically equivalent ways.
The introduction of these valuations gives an algebraic cast to the theory.
If one has a fact about matroids stated in elementary set-theoretic terms
and wishes to produce its counterpart for valuated matroids, 
the procedure is far from mechanical, 
but as a rule one reads the set-theoretic operations as boolean algebra
and then changes coefficients from the boolean semifield 
$\mathbb B=(\{\text{False},\text{True}\}, \text{OR}, \text{AND})$
to the \emph{tropical semifield}
$\mathbb T=(\mathbb R\cup\{\infty\}, \min, +)$,
where in these triples we list the addition operation second and the multiplication third.
The symbol $\infty$ denotes positive infinity, so $\min(\infty,a)=a$ and $\infty+a=\infty$ for all $a\in\mathbb T$.
Of these new coefficients $\infty$ behaves like False, 
and $a\in\mathbb R$ like True with $a$ attached as a valuation.
(To make a precise statement,
there is a semifield homomorphism $\mathbb T\to\mathbb B$
with $\infty\mapsto\text{False}$ and $a\mapsto\text{True}$ for $a\in\mathbb R$.)

As a first example, presentations of transversal valuated matroids 
are bipartite graphs with real edge weights.
The boolean predicate ``$(u,v)$ is an edge'' 
becomes a weight in $\mathbb T$ on the pair $(u,v)$, with $\infty$ encoding the case that the edge $(u,v)$ is absent.

The authors' previous work \cite{FinkOlarte} investigated 
the set of all presentations of a valuated matroid.
Our version of Brualdi and Dinolt's result above is cited as \Cref{FO6.6} below.
Where they have minimal presentations, 
we have presentations of minimal \emph{support}, i.e.\ where the graph is minimal but the real valuations are not otherwise constrained.
We are now ready to state two of our new main results.
This paper's \Cref{thm:different} generalizes \cite[Theorem 3.1]{BdM}
by showing that if a presentation is of minimal support, 
then adding a new vertex yields a different valuated matroid extension 
for each choice of neighbourhood and of real weights on the new edges.
\Cref{thm:minimal} generalizes \cite[Theorem 3.3]{BdM} by showing that every transversal valuated extension can be obtained by adding a vertex to a presentation of minimal support. 

To motivate the shape of our last main result, 
we explain a second example of translating from matroids to valuated matroids,
namely flats.
Among the flat axioms for matroids
is the property that the flats form a lattice under intersection,
containing the ground set $E$.
As finite semilattices are lattices, it is enough to ask just for a $\cap$-semilattice,
and we can fold in the condition that $E$ is a flat by 
requiring closure under intersection of arbitrary arity, including arity~0.
For the correct valuated generalization it turns out we need to think of the complements of flats, which are closed under arbitrary unions.
Set union is coordinatewise OR of boolean vectors, i.e.\ coordinatewise sum over $\mathbb B$.
When we change coefficients to~$\mathbb T$,
we should ask for a set of vectors in~$\mathbb T^E$ closed not just under addition 
but under all $\mathbb T$-linear combinations;
after all, an arbitrary $\mathbb B$-linear combination of some vectors is just the sum in~$\mathbb B$ of a subcollection.
Such a set of vectors in~$\mathbb T^E$ is called a $\mathbb T$-\emph{semimodule}, 
the semifield analog of a module over a ring.
We will not talk through the translation of the last flat axiom, the partition property for covers;
this becomes what is called the ``balancing condition'' in tropical geometry,
implying that our semimodule is a \emph{tropical linear space}.

It is known that the single-element extensions of a fixed matroid form a lattice under containment of the families of independent sets. Bonin and de Mier ask whether transversal extensions form a sublattice.
As above, the corresponding question in the valuated setting is whether transversal extensions of valuated matroids form a $\mathbb T$-semimodule.
Our \Cref{thm:join} shows that extensions arising from a single minimal presentation form a $\mathbb T$-semimodule,
generalizing \cite[Theorem 4.4]{BdM}. But the general question remains open. 

\subsection*{Acknowledgments} The first author received support from 
the Engineering and Physical Sciences Research Council [grant number EP/X001229/1] 
and the Deutsche For\-schungs\-ge\-mein\-schaft project ``Facetten der Komplexit\"at''.
The second author was supported by project CLaPPo (21.SI03.64658) of Universidad de Cantabria and Banco Santander.

\section{Definitions and notation}

We write multisets between double braces, like $\ldb0,0,1\rdb$.
The multiplicity with which an element $p$ appears in a multiset $A$
is denoted $\mult(p,A)$.

\subsection{Matroids}
See \cite{Oxley} for basic definitions in matroid theory, including those of the terms in this paragraph.
Given a matroid $M$ on the finite set $E=E(M)$,
we let $\BB(M)$ denote the set of \newword{bases} of~$M$,
$\CF(M)$ its lattice of \newword{cyclic flats},
$\rk(M)$ its \emph{rank},
and $\rk_M(S)$ and $\cork_M(S)$ the \newword{rank} and \newword{corank} of a set $S\subseteq E(M)$.

\subsection{Valuated matroids}
A \newword{valuated matroid} of \newword{rank}~$d$ is a function $\mu$ 
from the set $\binom Ed$ of $d$-element subsets of~$E=E(\mu)$ to $\mathbb T=\mathbb R\cup\{\infty\}$,
the value of $\mu$ at~$S\in\binom Ed$ being written $\mu_S$, such that
\begin{enumerate}
\item $\mu^{-1}(\mathbb R) = \BB(\underm\mu)$ for some matroid $\underm\mu$,
called the \emph{underlying matroid} of~$\mu$, and
\vspace{0.5ex} 
\item $\mu$ satisfies the \newword{3-term Pl\"ucker relations\/}: the minimum
\end{enumerate}
\[
\min(\mu_{Sij}+\mu_{Skl},\ \mu_{Sik}+\mu_{Sjl},\ \mu_{Sil}+\mu_{Sjk})
\]
is achieved by at least two of the arguments, for each $S\in \binom E{d-2}$ and $i,j,k,l \notin S$.
The notation e.g.\ $Sij$ is short for $S\cup\{i,j\}$. 

However, distinct such functions $\mu$ and $\nu$ 
may represent the same valuated matroid. This is the case
if and only if there exists $t\in\mathbb R$ such that $\nu_S = \mu_S+t$ for all $S\in\binom Ed$.
We will write a valuated matroid simply as $\mu$, rather than with extra notation for the implicit equivalence class of equivalent functions.
The precise choice of representative function $\mu:\binom Ed\to\mathbb T$ will be important especially in the later parts of \Cref{sec:minimal}.

If $F\subseteq E(\mu)$ is a set such that $E(\mu)\setminus F$ is spanning in $\underm{\mu}$, then
the \emph{deletion} of $F$ from the valuated matroid $\mu$
is the valuated matroid $\mu\setminus F$ obtained as the restriction of the function $\mu$ to $\binom{E\setminus F}d$.
Similarly, if $F$ is independent in $\underm{\mu}$, then 
the \emph{contraction} $\mu /F$ is the function $B \mapsto \mu(B\cup F)$ for every $B\in \binom{E\setminus F}d$.
The single-element deletion $\mu\setminus\{j\}$ and contraction $\mu/\{j\}$ are more often written $\mu\setminus j$ and $\mu/ j$, respectively.

Given a valuated matroid $\mu$, its \newword{tropical linear space} 
$\Be \mu \subseteq \mathbb T^{\nE}$ 
consists of all points $x$ where the minimum 
\[
\min_{i\in T} (\mu_{T\setminus i} + x_i)
\]
is achieved at least twice for every $T\in{E \choose d+1}$. When $\mu$ is the trivial valuated matroid over $M$, 
i.e.\ it is identically 0 on bases of~$M$, then $\Be \mu =: \Be M$ is called the \newword{Bergman fan} of $M$.

The (basis) \newword{polytope} $P_M$ of a matroid $M$ is the convex hull of the $\{0,1\}$ indicator vectors of its bases.
A valuated matroid $\mu$ induces a regular subdivision of the polytope of $\underm\mu$ into matroid polytopes \cite[Lemma 4.4.6]{MaclaganSturmfels}. 
We write $\MM(\mu)$ for the set of all matroids $M$ whose polytopes appear in this subdivision
such that every loop of~$M$ is a loop of~$\underm{\mu}$. 
All matroids $M\in\MM(\mu)$ satisfy $\rk(M)=\rk(\underm{\mu})=\rk(\mu)$
and $M\leq\underm{\mu}$ in the weak order
(that is, independent sets of $M$ are independent in $\underm{\mu}$).
For $x\in\Be \mu\cap\mathbb R^E$, the set of $S\in \binom Ed$ which minimize $\mu_S+\sum_{i\in S}x_i$
is the set of bases of a matroid $M\in\MM(\mu)$.
($M$ is known as the \newword{initial matroid} of $\mu$ at~$x$.)
The locus of points $x\in\Be \mu$ yielding a given $M$ is a relatively open polyhedron, whose closure we denote $\f {\Be \mu} M$.
As $M\in\MM(\mu)$ varies, the $\f {\Be \mu} M$ are the faces of a polyhedral complex structure on~$\Be \mu\cap\mathbb R^E$.
In a small neighborhood of any point $x\in \f {\Be \mu} M$, 
the tropical linear space $\Be \mu$ agrees with the translate $x + \Be M$.


The \newword{support} of a vector $y\in\mathbb T^E$ 
is $\supp(y) := \{j\in E:y_j\ne\infty\}$.
The support of a multiset $\AAA$ of such vectors is the multiset $\ldb\supp(p):p\in\AAA\rdb$.


\subsection{Extensions}

\begin{definition}\label{def:valuated extension}
A (single-element) \newword{extension} of a valuated matroid $\mu$ on~$\nnE$ of rank~$d$,
by an element $\nplusone$, is a
valuated matroid $\mu'$ on $\nnEplusone$ of rank~$d$ such that $\mu'\backslash \nplusone = \mu$.
\end{definition}
A nonvaluated matroid has a unique extension that increases its rank, namely the addition of a coloop;
all other extensions preserve the rank.
Similarly there is a unique $\mu'$ of rank $d+1$ such that $\mu'\backslash \nplusone = \mu$, in which $\nplusone$ is a coloop.
The remainder of this paper ignores that case. From now on, all extensions preserve rank.

\begin{remark}
Although we will not need them in this paper, we briefly recount the valuated analogues of other cryptomorphic interpretations of matroid extensions. See \cite[\S4.2.2]{Frenk} for a full treatment.

Single-element extensions of a matroid $M$ are in bijection with corank~1 quotients $M\to M'$.
In the valuated matroid context, this generalizes to the fact that extensions of $\mu$ are cryptomorphic to corank 1 tropical linear subspaces of $\Be \mu$. To see this, notice that the Pl\"ucker relations of an extension $\mu'$ are equivalent to the incidence relations between $\mu$ and $\mu'/\nplusone$.

Another object in bijection with single-element extensions of an unvaluated matroid $M$ are modular cuts. 
This generalizes in the valuated setting to \emph{extension functions} $f : \Be \mu \to \TT$,
which are the functions on $\Be\mu\subseteq\mathbb T^{\nE}$ linear in an appropriate sense;
in particular they include all restrictions of global linear functions on $\mathbb T^{\nE}$ (corresponding to principal modular cuts). 
Given an extension $\mu'$ of $\mu$, we have that $\Be {\mu'}$ maps to $\Be\mu$ 
under the projection $\pi : \TT^{\nEplusone} \to \TT^{\nE}$. The extension function associated to $\mu'$ is defined as $f(x) := \min \pi^{-1}(x) \cap \Be {\mu'}$. Knowing the values taken by $f$ on the cocircuits of $\mu'$ determines all of $f$.

The set of all extensions of a fixed matroid has a lattice structure (in the weak order).
For valuated matroids, the pointwise minimum of extension functions for~$\mu$ is again an extension function,
and adding a global scalar preserves the property of being an extension function.
That is, the family of extension functions for~$\mu$ is a $\mathbb T$-semimodule.
As in our introduction, this is the expected valuated generalization of a lattice.
\end{remark}

\section{Transversal (valuated) matroids and presentations}

Now we recall the key constructions related to presentations of transversal valuated matroids,
following \cite{FinkOlarte}. Before that, we briefly discuss the unvaluated case.

Given a bipartite graph $G$ with vertices $\nnE\sqcup [d]$, 
the subsets of $\nnE$ with a perfect matching to~$[d]$ in~$G$ are the bases of some matroid $M$ (unless there are no such subsets). 
Any matroid which can be constructed this way is called \newword{transversal}. 
In fact we will not work with graphs
but will encode $G$ by a multiset $\AAA= \ldb A_1 \dots A_d\rdb$ of subsets of~$\nnE$ where $A_i$ consists of the neighbours of $i\in[d]$;
this determines $G$ up to permutation of~$[d]$, and therefore $M$. 
Such a multiset $\AAA$ is called a \newword{presentation} of~$M$.
Our definition implies that the size of~$\AAA$ equals the rank of $M$,
and we maintain this assumption on presentations throughout the paper.

For any $A_i\in\AAA$, its complement $\nnE\setminus A_i$ is a flat in $M$:
a largest partial matching from $\nnE\setminus A_i$ to a subset of~$[d]$ cannot involve $i\in[d]$, and so for any element of~$j\in A_i$, the edge $(i,j)$ can be added to give a larger partial matching from $(\nnE\setminus A_i)\cup\{j\}$.

Given two presentations $\AAA$ and $\AAA'$ of $M$, we write $\AAA\le \AAA'$ if there is a bijection $f:\AAA\to \AAA'$ such that $X\subseteq f(X)$ for every $X\in \AAA$.
Under this order we have the following.
\begin{proposition}[{\cite[Theorem 4.2]{BrualdiDinolt}}]
\label{prop:BD}
Any transversal matroid $M$ has a unique maximal presentation. 
A presentation $\AAA$ of $M$ is minimal if and only if $A_i$ is a cocircuit of $M$ for every $i$.
\end{proposition}
Brualdi and Dinolt also provide an algorithm to obtain all presentations \cite[\S4]{BrualdiDinolt}.

We move on to the valuated case.
\begin{definition}\label{def:valuated presentation}
A \newword{presentation} $\AAA$ of a valuated matroid $\mu$ on~$\nnE$ is a multiset 
$\ldb A_1,\ldots,A_d\rdb$ of vectors $A_i\in\mathbb T^{\nE}$
such that, for all sets $B\in\binom{\nnE}d$,
\begin{equation}\label{eq:mu_B}
\mu_B=\min\left\{\sum_{i\in B} (A_{\phi(i)})_i : \phi \text{ is a bijection } B\to[d]\right\}.
\end{equation}
Any valuated matroid $\mu$ having a presentation is called \newword{transversal}.
\end{definition}
If $\AAA$ is a presentation of $\mu$, then the support of~$\AAA$ is a presentation of $\underm{\mu}$.
This is the fact that justifies calling matroids in this class ``transversal valuated matroids'':
their underlying matroids are transversal in the ordinary sense.

In \Cref{sec:minimal} we will speak of presentations in the language of tropical matrices
that appears in~\cite{FinkOlarte} (in turn taken from \cite{FR}).
The \newword{tropical Stiefel map} $\stiefel$ 
maps a $d\times\nE$ matrix $A$ with entries in $\mathbb T$, whose rows we call $A_1,\ldots,A_d$, 
to the function $\mu:\binom Ed\to\mathbb T$ given by \eqref{eq:mu_B},
which are viewed as tropical maximal minors of~$A$.
$\stiefel(A)$ is always a valuated matroid as long as it is not the constant function $\infty$ \cite[p.~301]{FR}. 
In this way we can think of the tropical Stiefel map $\stiefel$ as acting on a multiset of 
vectors in $\mathbb T^{\nE}$,
and the presentations of $\mu$ as the elements in its preimage under $\stiefel$. 
Every row of~$A$ is an element of~$\Be{\stiefel(A)}$ \cite[Lemma 4.1]{FinkOlarte}.


\begin{remark}\label{rmk:projective}
In \cite{FinkOlarte} the elements of presentations, and of tropical linear spaces, 
were instead taken to be points of a \emph{tropical projective space},
the quotient of $\mathbb T^{\nE}$ which identifies vectors $x$ and $y$ if there is $\lambda\in\mathbb R$ such that $y_i=x_i+\lambda$ for all $i\in\nnE$.
The main advantage of doing so is that the elements of $\dapx(\mu)$ below are truly uniquely determined by $\mu$,
rather than being determined up to an additive scalar.
This makes the geometry tidier, notably in \Cref{FO6.6}.
But for the purposes of the present paper the distinction is not important.
\end{remark}

If $M$ is a transversal matroid and $F$ a cyclic flat of~$M$, then $M/F$ is again transversal \cite[Theorem 5.4]{BrualdiMason}.
Likewise, if $\mu$ is a transversal valuated matroid and $F$ a cyclic flat of $\underm{\mu}$,
then $\mu/F$ is transversal \cite[Proposition 4.14]{FinkOlarte}.

\begin{definition}
Assume $\mu$ is transversal.
The \newword{distinguished multiset of matroids} $\dmat(\mu)$ is defined as follows.
For each cyclic flat $F$ of $\underm\mu$ and each $M\in\MM(\mu/F)$,
we have that $\dmat(\mu)$ contains $M$ with multiplicity
\[t(M):=\sum_{F'\in\CF(M)}\Moebius(\emptyset,F')\cork_M(F'),\]
where $\Moebius$ is the M\"obius function on $\CF(M)$.
\end{definition}

What we call $t(M)$ is denoted~$\tau_M(\emptyset)$ in~\cite{FinkOlarte}.
The integers $t(M)$ are non-negative
by \cite[Theorem 4.7]{BrualdiDinolt} and \cite[Corollary 5.6]{FR}, 
in view of the fact that $\mu/F$ is transversal.
A necessary condition for $t(M)>0$ is that $M$ is connected. 

\begin{proposition}[{\cite[Proposition 6.2]{FinkOlarte}}]
We have $|\dmat(\mu)|=d$.
\end{proposition}

Given $M\in\dmat(\mu)$, there is a unique cyclic flat $F(M)$ of~$\underm{\mu}$ such that $M\in \MM(\mu/F(M))$.
The cell $\f{\Be{\mu/F(M)}}M$ is one-dimensional.
Let $\vertex\mu M\in\mathbb R^{\nnE\setminus F(M)}$ be 
any vector in $\f{\Be{\mu/F(M)}}M$.
Then $\f{\Be{\mu/F(M)}}M=\{\vertex\mu M+(\lambda,\lambda,\ldots,\lambda):\lambda\in\mathbb R\}$.

\begin{definition}
    Given a transversal valuated matroid $\mu$, the \emph{distinguished multiset of apices} $\dapx(\mu)$ consists of $\overline{\vertex\mu M}$ for every $M\in\dmat(\mu)$, with the same multiplicity, where 
$\overline{\vertex\mu M}\in \mathbb T^{\nE}$ 
extends $\vertex\mu M$ by supplying $\infty$ for the missing coordinates (those indexed by~$F(M)$). 
\end{definition}

Let $e_j$ be the $j$\/th standard coordinate vector in $\mathbb R^n$,
with $j$\/th coordinate 1 and other coordinates 0.
We now present the main structural theorems from \cite{FinkOlarte} on the set of presentations of a valuated matroid.
\Cref{FO6.6} has been rewritten to avoid unintroduced notation.

\begin{theorem}\label{thm:FOmaximal}
If $\mu$ is a transversal valuated matroid, then $\dapx(\mu)$ is a presentation of $\mu$,
and its support is the maximal presentation of $\underm{\mu}$.
\end{theorem}

\begin{proof}
That $\dapx(\mu)$ is a presentation of~$\mu$ is a consequence of~\cite[Theorem 6.6]{FinkOlarte}.
To spell this out, the elements $\overline{\vertex\mu M}$ of $\dapx(\mu)$
are the same as the $\vertex LM$ appearing in \cite[Definition 6.3]{FinkOlarte},
and the tuple consisting of $t$ copies of the zero vector lies in the set $\phi(M)$ of that definition
(the maximal presentation of $M$ witnesses this containment, \cite[Theorem 4.7]{BrualdiDinolt}).

For a set $F\subseteq E(\mu)$, 
all vectors of $\dapx(\mu)$ of support $E(\mu)\setminus F$
are of the form $\overline{\vertex\mu M}$ for some $M\in\MM(\mu/F)$. 
The total multiplicity of all these $\overline{\vertex\mu M}$ is $t(\underm{\mu}/F)$
\cite[Proof of Proposition 6.2]{FinkOlarte}.
This equals the multiplicity with which $E(\mu)\setminus F$ appears in the maximal presentation of~$\underm{\mu}$
\cite[Theorems 3.5, 3.8]{BoninTransversalNotes}.
\end{proof}

\begin{theorem}[{\cite[Theorem 6.6]{FinkOlarte}}]\label{FO6.6}
Let $\mu$ be a transversal valuated matroid.
An arbitrary multiset $\AAA$ of vectors in~$\mathbb T^{\nE}$ is a presentation of~$\mu$ if and only if it is a disjoint union
\[\AAA=\coprod_{M\in\dmat(\mu)}\AAA_M\]
satisfying the following conditions:
\begin{enumerate}
\item $|\AAA_M|=\mult(M,\dmat(\mu))$;
\item For each $p\in \AAA_M$ we have $p= \overline{\vertex\mu M} +(\lambda,\lambda,\ldots,\lambda)+\sum\limits_{j\in J_p}\alpha_j e_j$ for some $J_p$ an independent flat in~$M$,
$\alpha_j\in [0,\infty]$, and $\lambda\in\mathbb R$;
\item Some presentation of~$M$ contains the submultiset $\ldb E(M)\setminus J_p : p\in \AAA_M\rdb$.
\end{enumerate}
\end{theorem}

To be precise about the rewriting,
\cite[Theorem 6.6]{FinkOlarte} invokes a set $\Pi(L)$, which is defined in \cite[Definition 6.3]{FinkOlarte}.
We have renamed the factor $\phi_L(M)$ in the definition of~$\Pi(L)$ to $\AAA_M$,
and because the present paper uses multisets rather than matrices with ordered rows we do not need a symmetric group action.
We have replaced $\phi(M)$ by its definition,
and the use of ``relative support'' $\mathrm{rs}_0$ \cite[Definition 4.2]{FinkOlarte} by 
the parametrization $\sum_{j\in J_p}\alpha_j e_j$ of the vectors with the required relative support.
(When $M$ is the free matroid on~$E(M)$, the independent flat $J_p=E(M)$ cannot arise from the definition of relative support;
our statement excludes the case $J_p=E(M)$ through condition (3).)

\Cref{thm:FOmaximal} suggests that $\dapx(\mu)$ plays the role of the unique maximal presentation of a usual valuated matroid.
Item (2) of \Cref{FO6.6} is another expression of ``maximality'' of $\dapx(\mu)$.
Just as other presentations of a transversal matroid are obtained by deleting edges from the maximal presentation,
other presentations of a transversal valuated matroid are obtained by increasing coordinates in~$\dapx(\mu)$, potentially replacing them with $\infty$ which deletes them from the support.
The maximal presentation of a transversal matroid is unique,
while $\dapx(\mu)$ is the unique multiset which works in \Cref{FO6.6} up to the indeterminacy in the choice of $\vertex\mu M$ (accounted for by \Cref{rmk:projective}).

The next lemma shows that the submultisets $\AAA_M$ in \Cref{FO6.6} are functions of~$\AAA$, 
so that the notation ``$\AAA_M$'' can be used unambiguously once $\AAA$ is chosen.
\begin{lemma}\label{lem:A_M well defined}
For any presentation $\AAA$ of~$\mu$,
there is a unique decomposition $\AAA=\coprod_M \AAA_M$ meeting the conditions of \Cref{FO6.6}. 
\end{lemma}

\begin{proof}
This is implicit in the proof of \cite[Proposition 6.7]{FinkOlarte}.
Let $F$ be a cyclic flat of~$\underm{\mu}$.
Let $\AAA_F$ be the submultiset of $\AAA$ of vectors whose support is disjoint from~$F$.
Because $(\overline{\vertex\mu M})_j=\infty$ for $j\in F$,
\Cref{FO6.6}(2) implies that $\AAA_M$ must be a submultiset of $\AAA_{F(M)}$ for all $M\in\dmat(\mu)$.

By \cite[Proposition 4.14]{FinkOlarte},
the projection of $\AAA_F$ to the coordinates indexed by $\nnE\setminus F$
gives a presentation of $\mu/F$. Call the projection $\AAA/F$.
Let $M\in\dmat(\mu)$ satisfy $F(M)=F$, and let $G$ be a cyclic flat of~$M$.
By \cite[Proposition 4.3(2)]{FinkOlarte} applied to 
the ``zoomed'' valuated matroid denoted $Z_{\vertex\mu M}(\AAA/F)$ there,
the number of vectors in $\AAA/F$ of the form 
\begin{equation}\label{eq:FO6.7}
\vertex\mu M +(\lambda,\lambda,\ldots,\lambda)+\sum\limits_{j\in J\setminus G}\alpha_j e_j
\end{equation}
for $\lambda\in\mathbb R$, $J$ a set containing $G$, and each $\alpha_j>0$ strictly,
is $\cork_M(G)$.
By M\"obius inversion, there are exactly $t(M)=\mult(M,\dmat(\mu))$ vectors in $\AAA/F$ of the form \eqref{eq:FO6.7} with $J$ independent.
By \Cref{FO6.6}(1), the image of $\AAA_M$ in $\AAA/F$ must be of size $t(M)$,
so by \Cref{FO6.6}(2) it must be exactly the multiset of the $t(M)$ vectors in the previous sentence.
This determines $\AAA_M$.
\end{proof}

\begin{definition}\label{def:minimal}
A presentation $\AAA$ of $\mu$ is called \newword{minimal} if the support of $\AAA$ is a minimal presentation of $\underm{\mu}$.
\end{definition}

We use the following proposition.

\begin{proposition}\label{prop:minimal}
A presentation $\AAA$ of $\mu$ is minimal if and only if for every $M\in\dmat(\mu)$ and every $p\in\AAA_M$,
the set $E(M)\setminus\supp p$ is a hyperplane of~$M$.
\end{proposition}

\begin{proof}
For $p\in\AAA_M$ let $F(p)$ be the set of infinite coordinates of $p$. 
By \Cref{FO6.6}(2), we have $F(M)\subseteq F(p)$.

By \Cref{prop:BD}, $\AAA$ is minimal if and only if $F(p)$ is a hyperplane of $\underm{\mu}$ for each $p$.
Since $F(p)\setminus F(M)$ must be an independent flat of~$M$ by \Cref{FO6.6}(2), we have
$\rk_{M}(F(p)\setminus F(M)) = |F(p)\setminus F(M)|$.
We have that $\rk(M) = \rk(\underm{\mu})- \rk_{\underm{\mu}}(F(M))$
and $M\leq \underm{\mu}/F(M)$ in the weak order.
The latter implies that $F(p)$ is also independent of $F(M)$ in~$\underm{\mu}$,
giving the second equality in
\begin{align*}
\rk(M) - \rk_{M}(F(p)\setminus F(M)) &= \rk(\underm{\mu})- \rk_{\underm{\mu}}(F(M)) - |F(p)\setminus F(M)| 
\\&=\rk(\underm{\mu})-\rk_{\underm{\mu}}(F(p)).    
\end{align*}
So $F(p)$ is a hyperplane of~$\underm{\mu}$ if and only if $F(p)\setminus F(M)$ is a hyperplane of~$M$.
By definition $F(p)\setminus F(M)=E(M)\setminus\supp p$.
\end{proof}

\section{Minimal transversal presentations and their extensions}\label{sec:minimal}

We study single element extensions of transversal valuated matroids which are also transversal. 
In this section it will be convenient to encode all transversal valuated presentations $\AAA$
as tropical matrices $A$, which we'll call \emph{matrix presentations},
whose rows are the elements of~$\AAA$.
We call a matrix presentation minimal if the equivalent multiset presentation is.
If $A\in\mathbb T^{d\times\nE}$ is such a matrix,
we write $(A|x)\in\mathbb T^{d\times(\nE\cup\{\nplusone\})}$ for the matrix with rows $(A_i,x_i)$, the last column being indexed by a new element $\nplusone$.
Recall the map $\stiefel$ from the discussion after \Cref{def:valuated presentation}.

\begin{definition}\label{def:valuated transversal extension}
	Let $\mu$ be a transversal valuated matroid.
    A \newword{transversal extension} of $\mu$ is a valuated matroid of the form $\mu'= \stiefel(A|x)$ for a matrix presentation $A$ of $\mu$ and a vector $x\in \TT^d$. 
\end{definition}


If $A$ is a matrix presentation of a valuated matroid $\mu'$, then removing the $j$\/th column 
of~$A$ yields a matrix presentation of $\mu'\backslash j$. 
Therefore for any transversal extension $\mu'$ of a valuated matroid $\mu$ on~$\nnE$, 
there is some presentation 
$\AAA=\ldb A_1,\ldots,A_d\rdb$ of~$\mu$
so that concatenating a new tropical number to each $A_i$ produces a presentation of $\mu'$.
In other words, any extension (\Cref{def:valuated extension}) of a transversal valuated matroid which is again transversal
is indeed a transversal extension in the sense of \Cref{def:valuated transversal extension}.

However, one cannot obtain every transversal extension $\mu'$ by extending a fixed matrix presentation $A$ of $\mu$, as the following example shows.
\begin{example}
\label{ex:u23}
Consider the rank 2 matroid $\mu$ on 3 elements where $\mu(B) = 0$ for every $B\in {[3] \choose 2}$.
For $i\in[3]$ and $\lambda>0$, the valuated matroid $\mu^{i,\lambda}$ on $[4]$ which is 0 everywhere except that $\mu^{i,\lambda}_{i4}=\lambda$ is a transversal extension of $\mu$.
Because $\mu^{i,\lambda}_{i4}>\mu^{i,\lambda}_{j4}$ for $j\in[3]\setminus\{i\}$, in any matrix presentation of $\mu^{i,\lambda}$ there is a row where the $i$th entry is larger than the other entries in the first three columns.
To do a particular case concretely, $\mu^{1,1}=\stiefel(A|x)$ implies that, up to exchanging the two rows and adding constants to rows or columns,
\[
A|x=\left(\begin{array}{@{\,}ccc|c@{\,}}a&0&0&c\\0&b&0&0\end{array}\right) 
\quad\mbox{or}\quad
A|x=\left(\begin{array}{@{\,}ccc|c@{\,}}a&0&0&c\\0&0&b&0\end{array}\right)
\]
for some $a,b,c\ge0$ with $\min(a,c)=1$, in particular $a\ge1$.
But as any fixed matrix presentation $A$ of $\mu$ has only two rows, 
there are only two $i$ such that $\mu^{i,\lambda}= \stiefel(A|x)$ is possible.

Note that, by taking $\lambda=\infty$, 
the example shows that not even all underlying matroids of transversal extensions of~$\mu$ 
can be obtained by extending a fixed presentation of~$\mu$.
\end{example}

Given a transversal extension of~$\mu$, 
it is then a natural question from which presentations $A$ of~$\mu$ it can be obtained.
Our answers are parallel to those of Bonin and de Mier for the unvaluated case.
We begin with the valuated counterpart of \cite[Theorem 3.1]{BdM}.
\begin{theorem}
\label{thm:different}
Let $A$ be a matrix presentation of a valuated matroid $\mu$. 
Then all extensions $\stiefel(A|x)$ for $x\in \TT^d$ are different if and only if $A$ is minimal.
\end{theorem}

\begin{proof}
By \Cref{lem:A_M well defined}, the distinguished matroids of~$\mu$ may be named $M_1,\dots,M_d$, with the correct multiplicities,
such that the vector $A_i$ lies in $\AAA_{M_i}$ for all $i\in[d]$.

Suppose $A$ is not minimal. 
By \Cref{prop:minimal} there exists $i\in[d]$ such that 
$E(M_i)\setminus\supp A_i$ is not a hyperplane of~$M_i$. 
Consider $x_t= te_i\in \TT^d$, the vector whose $i$-th coordinate is a real number $t$ and whose other coordinates are~$0$.
Then each entry of $\stiefel(A|x_t)$ has the form $\min\{a+t,b\}$ for $a,b\in\mathbb T$
(i.e.\ is a tropical affine function),
because the weight of each matching is either a constant or a constant plus~$t$.
We are going to show that the constant term $b$ never equals $\infty$ unless also $a=\infty$.
So $\stiefel(A|x)$ stabilizes for a sufficiently large $t$.

Suppose that instead there is a set $J\subseteq E$, with $|J|=d-1$, such that $\stiefel(A|x_t)_{J\cup\{\nplusone\}} = a+t$ for some $a\in\mathbb R$. So $J$ is independent in $\underm{\mu}$ and $A$ supports a matching $\phi: J\cup\{\nplusone\}\rightarrow [d]$ with $\phi(\nplusone) = i$, i.e.\ $A_{\phi(j),j}$ is finite for $j\in J$. If there is any $j\in J$ such that $A_{i,j}$ is finite, then swapping the values of $j$ and $\nplusone$ in $\phi$ yields a matching on finite entries not using $t$, producing a constant term in the minor $(A|x)_{J\cup\{\nplusone\}}$. So $A_{i,j} = \infty$ for all $j\in J$. As $J$ is independent in $\underm{\mu}$, it follows that 
$$
|J\backslash E(M_i)| = \rk_{\underm{\mu}}(J\backslash E(M_i)) \le \rk_{\underm{\mu}}(\nnE\backslash E(M_i)) = d-\rk(M_i),
$$
where the last equality follows from $\rk(M_i) = \rk(\underm{\mu}/(\nnE\backslash E(M_i))$. Therefore
$$
|J\cap E(M_i)| = |J|-|J\backslash E(M_i)| \ge d-1 -(d-\rk(M_i)) = \rk(M_i)-1.
$$
As $\supp(A_i)\cap J = \emptyset$, we have that $E(M_i)\backslash \supp(A_i)$ contains $J\cap E(M_i)$. But $E(M_i)\backslash \supp(A_i)$ is independent in $M_i$ of rank at most $\rk(M_i)-1$
(because it is a flat and $\supp(A_i)\ne\emptyset$).
So it must also be a hyperplane, which contradicts our initial assumption.

Now for the other direction let $A$ be a minimal matrix presentation of $\mu$. Let $\mu(x) = \stiefel(A|x)$. 
Take any $B_0\in\binom{\nnE}d$ such that $\mu_{B_0}$ is finite; then $\mu(x)_{B_0}=\mu_{B_0}$ for any~$x$.
We are going to show that for every $i\in[d]$ there is a maximal minor $(A|x)_{B_i}$, for $B_i \in \binom{E \cup \{*\}}{d}$,
that is of the form $a_i+x_i$ for some $a_i\in\mathbb R$.
If we do this, then for $x,y\in\mathbb T^d$ we have that
$x_i\ne y_i$ implies that $\mu(x)_{B_i}-\mu(x)_{B_0}\ne \mu(y)_{B_i}-\mu(y)_{B_0}$ and hence $\mu(x) \ne \mu(y)$. 
Therefor all extensions $(A|x)$ differ for different $x$. 

Consider $H_i = E(M_i)\backslash \supp(A_i)$, which is a hyperplane of $M_i$ by \Cref{prop:minimal}
and independent by \Cref{FO6.6}(2).
Let $F_i=\nnE\backslash E(M_i)$, which is a cyclic flat of $\underm{\mu}$. Let $J_i$ be a basis of the restriction of $\underm{\mu}$ to~$F_i$, and $B_i = J_i\cup H_i\cup\{\nplusone\}$.

As $|H_i|=\rk(M_i)-1=\rk(\underm{\mu}/F_i)-1=d-\rk_{\underm{\mu}}(F_i)-1$, we have that $|B_i|=d$.
Since $H_i$ is independent in $M_i$, it is independent in 
$\underm\mu/F_i$,
which has more independent sets because $M_i$ is among its initial matroids. 
So $B_i\backslash\{\nplusone\}$ is independent in $\underm{\mu}|_{F_i}\bigoplus \underm{\mu}/F_i$ and hence independent in $\underm{\mu}$. 
But all the entries $A_{i,j}$ with $j\in B_i\backslash\{\nplusone\}$ are infinite. 
Since $B_i\backslash\{\nplusone\}$ can be extended to a basis $B$ of $\underm{\mu}$, 
a matching $B\to[d]$ must exist,
and must restrict to $B_i$ as a matching $B_i\backslash\{\nplusone\}\to [d]\backslash\{i\}$. 
Hence there exists a matching of $\phi: B_i\to [d]$, but all such satisfy that $\phi(\nplusone)=i$.
This implies that $\mu(x)_{B_i} = \mu_{B_i\backslash\{*\}}+x_i$, as we wanted.  
\end{proof}

For completeness we record the counterpart of \cite[Lemma 3.2]{BdM}.
\begin{proposition}
Any extension $A|x$ of a minimal matrix presentation $A$ of a valuated matroid $\mu$
is itself a minimal matrix presentation of $\stiefel(A|x)$.
\end{proposition}

\begin{proof}
By \Cref{def:minimal}, this reduces immediately to \cite[Lemma 3.2]{BdM}.
\end{proof}

We state the next theorem of Bonin and de Mier that we will generalize,
with a fact that we will need from the proof extracted into the statement.
\begin{theorem}[{\cite[Theorem 3.3]{BdM}}]\label{thm:BdM3.3}
If $M$ is a transversal matroid on $\nnEplusone$ and $\nplusone$ is not a coloop,
then $M$ has a presentation that is an extension of some minimal presentation of $M\backslash\nplusone$.

In particular, if $\nplusone$ appears in exactly $k$ sets in the maximal presentation of~$M$,
then every minimal presentation of~$M$ in which $\nplusone$ appears in $k$ sets is an extension of a minimal presentation.
\end{theorem}

Here is our generalization.
\begin{theorem}
\label{thm:minimal}
Let $\mu'$ be a valuated transversal matroid on $\nnEplusone$ and suppose $\nplusone$ is not a coloop of $\underm{\mu'}$. Then $\mu'$ has a matrix presentation 
of the form $(A|x)$ such that $A$ is a
minimal presentation of $\mu' \backslash\nplusone$.
\end{theorem}
In other words, every transversal extension of a valuated matroid $\mu=\mu' \backslash\nplusone$ can be presented by adding a column to a minimal matrix presentation of~$\mu$.

\begin{proof}
By \Cref{thm:FOmaximal}, the support of $\dapx(\mu')$ is the maximal presentation of~$\underm{\mu'}$.
Let $k$ be the number of finite entries in the last column of the matrix of~$\dapx(\mu')$.
We show constructively that there is a minimal matrix presentation $A'$ of~$\mu'$ with $k$ finite entries in the last column.

For each $M\in\dmat(\mu)$, by \Cref{lem:A_M well defined} there is a well-defined submultiset $\dapx(\mu')_M$.
Label its elements $\dapx(\mu')_M=\ldb p(M)_1,\ldots,p(M)_{t(M)}\rdb$.
\Cref{FO6.6}(3) says that $\supp(\dapx(\mu')_M)$
is contained in a presentation $\AAA(M) = \ldb A(M)_1, \ldots, A(M)_{\rk m}\rdb$ of~$M$; order the elements so that 
\[\supp(\dapx(\mu')_M)=\ldb A(M)_1, \ldots, A(M)_{t(M)}\rdb\]
with $\supp(p(M)_i) = A(M)_i$. Applying \cite[Lemma 2.10]{BdM} to $\AAA(M)$,
there exists a minimal presentation $\AAA'(M) = \ldb A'(M)_1, \ldots, A'(M)_{\rk m}\rdb$ of~$M$
with $A'(M)_i\subseteq A(M)_i$ for all~$i$, and 
$\nplusone\in A'(M)_i$ if and only if $\nplusone\in A(M)_i$.

Define a matrix $A'\in\mathbb T^{d\times\nnEplusone}$ with a block of $t(M)$ rows for each $M\in\dmat(\mu')$: 
in the block for~$M$ are the vectors $p'(M)_1,\ldots,p'(M)_{t(M)}$ with
\[
(p'(M)_i)_j = \begin{cases}
(p(M)_i)_j & j\in A'(M)_i \\
\infty & \mbox{otherwise.}
\end{cases}
\]
Then $\supp(p'(M)_i) = A'(M)_i$
so, by construction, $A'$ has $k$ finite entries in the last column.
The set of positions $E(M)\setminus A'(M)_i$ of new $\infty$ coordinates in $p'(M)_i$
is an independent flat in~$M$ \cite[Theorem 4.7]{BrualdiDinolt},
so by \Cref{FO6.6}, $A'$ is a presentation of~$\mu'$.
By \Cref{prop:minimal}, $A'$ is a minimal presentation.


Write $A'=(A|x)$. Applying \Cref{thm:BdM3.3} to $\underm{\mu'}$, 
the support of~$A$ is a minimal presentation of $\underm{\mu'} \backslash\nplusone = \underm{\mu' \backslash\nplusone}$.
Therefore $A$ is a minimal presentation of $\mu' \backslash\nplusone$.
\end{proof}

\begin{definition}
    Given a transversal valuated matroid $\mu$ of rank $d$, 
    define the \emph{poset of transversal extensions} of~$\mu$ to have as its elements all transversal extensions of $\mu$ to~$\nnEplusone$.
		The order is given by $\mu'\le \mu''$ if $\mu'_B\le \mu''_B$ for every $B\in {\nnEplusone\choose d}$, under a choice of representatives of $\mu'$ and $\mu''$ that fixes $\mu'_B= \mu''_B$ for $B$ such that $\nplusone\notin B$.
\end{definition}
Readers used to the weak order on matroids beware: 
if $\mu'\le\mu''$ then $\underm{\mu'}\ge\underm{\mu''}$ in the weak order.
Because non-bases have $\infty$ as their valuation,
$\mu'_B\le \mu''_B$ translates to $B$ being a basis of $\mu'$ if it is a basis of $\mu''$. 

The following generalizes \cite[Theorem 4.4]{BdM}.
\begin{theorem}
\label{thm:join}
Let $\mu=\stiefel(A)$ be a valuated matroid and let $x, y \in \TT^d$.
Then in the poset of transversal extensions $\mu$,
$\stiefel(A|x)$ and $\stiefel(A|y)$ have a unique greatest lower bound, given by
\[
\stiefel(A|x) \vee \stiefel(A|y) = \stiefel(A|\min(x,y))
\]
where $\min(x,y)$ is meant coordinatewise.
\end{theorem}
\begin{proof}
Fix the matrix presentation $A$ of $\mu$. 
Given $J\in {\nnE\choose d-1}$, we have
\begin{equation}\label{eq:StAxJx}    
\stiefel(A|x)_{J\cup\nplusone} = \min_{j\in J} (A_{[d]\setminus j, J} + x_j)
\end{equation}
where $A_{I,J}$ is the tropical minor of~$A$ on row set~$I$ and column set~$J$. 
Any tranversal extension $\mu'$ of $\mu$ that is a common lower bound of $\stiefel(A|x)$ and~$\stiefel(A|y)$, 
must be coordinatewise bounded above by their coordinatewise minimum,
(where a representative of the coordinates of $\mu'$ is chosen so that $\mu'_B=\mu_B$ for $B$ such that $\nplusone\not\in B$).
Explicitly, this is
\[B\mapsto\begin{cases}
    \mu_B & \nplusone\not\in B \\
    \min\left(\min_{j\in J} (A_{[d]\setminus j, J} + x_j),
    \min_{j\in J} (A_{[d]\setminus j, J} + y_j)\right) \\
    \quad=\min_{j\in J} (A_{[d]\setminus j, J} + \min(x_j,y_j)) & \nplusone\in B,J=B\setminus\nplusone.
\end{cases}\]
Applying \Cref{eq:StAxJx} to $\min(x,y)$ shows that this function is the transversal extension $\stiefel(A|\min(x,y))$,
which means this is the unique greatest lower bound.
%
\end{proof}

The 
argument of the proof can be restated as follows.
Suppose that $\mu'$ and $\mu''$ are transversal extensions of $\mu$ to $\nnEplusone$
with respective matrix presentations $A', A''$ that agree on deletion of the last column,
representatives being chosen so that $\mu'_B= \mu''_B$ for $\nplusone\notin B$. Then
\[\mu'\vee\mu'' = \min(\mu',\mu'')\]
coordinatewise.
That is, for this choice of representatives, the meet in the poset of transversal extensions
becomes the tropical sum in the free $\mathbb T$-semimodule $\mathbb T^{\binom{\nnEplusone}d}$.
By the definition of valuated matroid, the set of all representatives of transversal extensions is also closed under tropical scalar multiplication.
But \Cref{thm:join} does not give a formula for the join in the case where $\mu'$ and $\mu''$ spring from different presentations of $\mu$.
\begin{question}
Let $\mu$ be a transversal valuated matroid on $\nnE$ of rank~$d$.
Is the set of all representatives in $\mathbb T^{\binom{\nnEplusone}d}$ of transversal extensions of $\mu$ to $\nnEplusone$ a sub-$\mathbb T$-semimodule of $\mathbb T^{\binom{\nnEplusone}d}$, 
with coordinatewise minimum as addition?
\end{question}

\bibliography{Stiefel_presentations2}

\begin{thebibliography}{10}

\bibitem{Bondy}
J.~A. Bondy.
\newblock Presentations of transversal matroids.
\newblock {\em J. London Math. Soc. (2)}, 5:289--292, 1972.

\bibitem{BoninTransversalNotes}
Joseph~E Bonin.
\newblock An introduction to transversal matroids, 2010.

\bibitem{BdM}
Joseph~E Bonin and Anna de~Mier.
\newblock Extensions and presentations of transversal matroids.
\newblock {\em European Journal of Combinatorics}, 50:18--29, 2015.

\bibitem{BrualdiMason}
Richard Brualdi and John Mason.
\newblock Transversal matroids and {H}all's theorem.
\newblock {\em Pacific Journal of Mathematics}, 41(3):601--613, 1972.

\bibitem{BrualdiDinolt}
Richard~A Brualdi and George~W Dinolt.
\newblock Characterizations of transversal matroids and their presentations.
\newblock {\em Journal of Combinatorial Theory, Series B}, 12(3):268--286,
  1972.

\bibitem{Crapo}
Henry~H. Crapo.
\newblock Single-element extensions of matroids.
\newblock {\em J. Res. Nat. Bur. Standards Sect. B}, 69B:55--65, 1965.

\bibitem{DressWenzel}
Andreas~WM Dress and Walter Wenzel.
\newblock Valuated matroids.
\newblock {\em Advances in Mathematics}, 93(2):214--250, 1992.

\bibitem{FinkOlarte}
Alex Fink and Jorge~Alberto Olarte.
\newblock Presentations of transversal valuated matroids.
\newblock {\em Journal of the London Mathematical Society}, 105(1):24--62,
  2022.

\bibitem{FR}
Alex Fink and Felipe Rinc{\'o}n.
\newblock Stiefel tropical linear spaces.
\newblock {\em Journal of Combinatorial Theory, Series A}, 135:291--331, 2015.

\bibitem{Frenk}
Bart Frenk.
\newblock {\em Tropical varieties, maps and gossip}.
\newblock PhD thesis, Technische Universiteit Eindhoven, 2013.

\bibitem{MaclaganSturmfels}
Diane Maclagan and Bernd Sturmfels.
\newblock {\em Introduction to tropical geometry}, volume 161.
\newblock American Mathematical Soc., 2015.

\bibitem{Mason-thesis}
John~Healey Mason.
\newblock {\em Representations of independence spaces}.
\newblock Phd thesis, The University of Wisconsin-Madison, 1970.

\bibitem{Oxley}
James Oxley.
\newblock {\em Matroid theory}, volume~21 of {\em Oxford Graduate Texts in
  Mathematics}.
\newblock Oxford University Press, Oxford, second edition, 2011.

\end{thebibliography}
\bibliographystyle{plain}

\end{document}